\documentstyle{article}
\title{Effective nonvanishing,\\ effective global generation }
\author{Mark Andrea A.  de Cataldo\thanks{
Partially supported by N.S.F. Grant DMS 9701779.}}

\newtheorem{tm}{Theorem}[section]
\newtheorem{lm}[tm]{Lemma}
\newtheorem{pr}[tm]{Proposition}
\newtheorem{rmk}[tm]{Remark}
\newtheorem{cor}[tm]{Corollary}
\newtheorem{ex}[tm]{Example}
\newtheorem{defi}[tm]{Definition}
\newtheorem{fact}[tm]{Fact}

\newtheorem{??}[tm]{Question}

\font\tenmsb=msbm10
\font\sevenmsb=msbm7
\font\fivemsb=msbm5

\newfam\msbfam
\textfont\msbfam=\tenmsb
\scriptfont\msbfam=\sevenmsb
\scriptscriptfont\msbfam=\fivemsb
\def\Bbb#1{{\fam\msbfam #1}}

\font\teneufm=eufm10
\font\seveneufm=eufm7
\font\fiveeufm=eufm5
\newfam\eufmfam
\textfont\eufmfam=\teneufm
\scriptfont\eufmfam=\seveneufm
\scriptscriptfont\eufmfam=\fiveeufm
\def\frak#1{{\fam\eufmfam\relax#1}}

\newcommand\ci{\cite}
\newcommand\s{\sigma}
\newcommand\rat{{\Bbb Q}}
\newcommand\comp{{\Bbb C}}
\newcommand\real{{\Bbb R}}

\newcommand\pn[1]{{\Bbb P}^{#1}}

\newcommand\blacksquare{{\hspace*{\fill} $\Box$}} 
\newcommand\odix[1]{ {\cal O}_{#1} }
\newcommand\odixl[2]{ {\cal O}_{#1}({#2}) }

\newcommand\e{\epsilon}

\newcommand\surj{-\!\!\!-\!\!\!-\!\!\!-\!\!\!\!\!\gg}
\newcommand\QQ{\frak Q}
\newcommand\q{\frak q}

\begin{document}

\maketitle

\begin{abstract}
We prove a multiple-points higher-jets  nonvanishing theorem by the use of
local Seshadri constants.  Applications are given to effectivity problems
such as constructing rational and birational maps into Grassmannians,
 and
the global generation of vector bundles.
\end{abstract}

\section{Introduction}
\label{intr}
Koll\'ar's nonvanishing 
theorem \ci{koebpf}, 3.2 is an instrument to make
Kawamata-Shokurov base-point-freeness assertion into an  effective one.
His result can be applied to a variety of other situations;
see \ci{koebpf}, \S4, \ci{koshafinv}, \S8 and \ci{koshaf}, \S14.

The basic set-up is as follows. Let $g: X \to S$ be a surjective
morphism of proper varieties, where $X$ in nonsingular and complete, $M$
be a  nef and $g$-big line bundle on $X$, $L$  be a nef and big line 
bundle on $S$
and $N=K_X + M + mg^*L$ be a line bundle  varying with the positive 
integer  $m$.
Koll\'ar proves, under the necessary assumption that
the first direct image sheaf  $g_*N\not= 0$,  that $h^0(X, N)=h^0(S,g_*N)>0$
and the sections 
of $g_*N$ generate this sheaf at a general point of $S$ for {\em every}
$\, m> (1/2)(\dim{S}^2 + \dim{S})$ (this is what makes the result into an 
``effective" tool).

The proof  starts with the choice of  a very general point $x$ on $S$ and 
ends 
with producing sections of $g_*N$ which generate at $x$ and therefore 
at a general point.

\medskip
The purpose of this note is to observe that we can  obtain
more precise statements  by considering the local Seshadri constants
of $L$ on $S$ and we can also simplify considerably the proof. See
the discussion at the beginning of  
\S\ref{nonvan} and Remark \ref{why}.
Our main result is the Effective Nonvanishing Theorem \ref{effnonvan},  
a ``multiple-points higher-jets"
version of \ci{koebpf}, Theorem 3.2.

The proof hinges on  Demailly's observation
 that given a nef line bundle $\cal L$ on $X$, a big enough local  
Seshadri constant 
 for $\cal L$ at a point $x$ can be used together with
Kawamata-Viehweg Vanishing Theorem to produce sections of the adjoint 
line bundle
$K_X + {\cal L}$ with nice 
generating properties  at $x$
(cf. \ci{dem94}, Proposition 7.10). An effective way to
force a big enough local Seshadri constant is Theorem \ref{ekltm}, which 
is due 
Ein-K\"uchle-Lazarsfeld. 

\smallskip
As applications we offer some generalizations (to the case of nef vector 
bundles)
 of the 
results  concerning line bundles in \ci{ekl} 
and  \ci{koshafinv}: effective construction of rational
and birational maps, and  nonvanishing on varieties with big enough 
algebraic 
fundamental group. These results in the case of one 
point  follow easily from results of Koll\'ar by considering the
tautological line bundle of the  projectivization of the vector bundles 
in question.
We present them as an exemplification of the unifying character
of Theroem \ref{effnonvan} and also because they are new in the case of
multiple-points and higher-jets.

\noindent
We also show how the global generation results for line bundles of
Anghern-Siu, Demailly, Tsuji and Siu (see \ci{dem94} and
\ci{deshm} for a  bibliography)
generalize to  vector bundles
of the form
$K_X^{\otimes a} \otimes E\otimes \det{E} \otimes L^{\otimes m}$, where
$a$ and $m$ are appropriate positive integers, $E$ is a nef vector bundle
 and $L$ is  and ample line bundle. 
We give explicit upper bounds on $m$ which   depend only on the dimension 
of the variety,
and not also on the Chern classes of the variety and the bundles in question.
However, we do not expect these bounds to be otpimal since 
they  do not match with the line bundle case ({\em i.e.} assuming that 
the vector bundle
$E$ is 
the trivial line bundle).

In our paper \ci{deshm} we prove upper  bounds as above for vector bundles
$E$ subject to
curvature conditions   which seem to be the natural  differential-geometric
analogue of nefness and indeed imply nefness.  These bounds match exactly
the results in the line bundle case. The methods employed there are 
analytical. 

In the final section we prove, using the language and techniques of algebraic Nadel ideals,  a global generating statement for nef vector bundles
 which indeed matches the result of Anghern and Siu in the line bundle
case.

\medskip
The paper is organized as follows.
\S1 is preliminary and consists of  easy and mostly known 
facts about local Seshadri constants, and of more elaborate ones, such as 
Theorem \ref{ekltm}, which makes 
Theorem \ref{effnonvan} into  an effective statement. 
\S2 is devoted to the main result, Theorem \ref{effnonvan}.
\S3 is devoted to the applications discussed above.
\S4 is devoted to the proof of Corollary 4.6,  a major improvment
of the results of Theorem 3.7.1 and 3.7.3.

\bigskip
\S4 has been written while the author enjoyed the hospitality
of the Max-Planck-Institut f\'ur Mathematik of Bonn.

\section{Notation and preliminaries}
\label{not}
We mostly employ the notation of
\ci{k-m-m}.

A {\em variety} is an integral separated scheme of finite type over
an uncountable algebraically closed   field of characteristic zero $k$.

We say that {\em a property holds at a very general point
on $X$} if it holds for every point  in the intersection, $\frak U$,
of some at most countable family of Zariski-open dense subsets of $X$. 
Any such set $\frak U$ meets any Zariski-open  dense subset of $X$.

The term ``point"  refers to   a closed one.

Vector bundles and associated locally free sheaves are identified. 
Cartier divisors are at times identified with the associated invertible 
sheaves
and 
the additive  and multiplicative
notation  are both used, at times simultaneously.

The symbol $B(a,b)$ denotes the usual binomial coefficient.

Let $X$ be a  
 variety, $n$ be its dimension and $Div(X)$  be the group
of Cartier divisors on $X$.
A {\em $\rat$-Cartier divisor} is an element of $Div(X)\otimes {\rat}$.
The linear  and numerical equivalence of $\rat$-divisors
 are  denoted by ``$\approx$"
and ``$\equiv$," respectively.  A {\em $\rat$-divisor} is an element
in $Z_{n-1}(X) \otimes {\rat}$, where $Z_{n-1}(X)$ is the free group
of Weil divisors on $X$.

The symbols $\lfloor a\rfloor$ and $\langle a \rangle$
denote the biggest integer less than or equal to $a$, 
and $a- \lfloor a\rfloor$, respectively. 
These symbols are used in conjunction with $\rat$-divisors when these 
divisors
are written as a $\rat$-combination of {\em distinct} prime divisors.

Given any proper morphism of varieties $\pi: X \to S$, we have the
 notions of ($\pi$-)ample, ($\pi$-)nef,  ($\pi$-)big and
($\pi$-)nef and ($\pi$-)big for (numerical equivalence
classes of)  $\rat$-Cartier divisors on $X$.

Let ${\cal A} \in (Div(X) \otimes {\rat})/\equiv$ be  a numerical  class 
 and $D$ be a $\rat$-Cartier divisor on $X$. 
By abuse of notation,
 ${\cal A} \equiv D$ means that $A \equiv D$
for one, and thus all, the elements in ${\cal A}$.
This remark plays a role when we use the canonical divisor class
 together with  $\rat$-Cartier divisors.

In what above, the field $\real$ can replace $\rat$ with minor changes.

\medskip
The following two vanishing-injectivity theorems are needed for 
Theorem \ref{effnonvan}.
\begin{tm}
\label{van}
{\rm (Cf. \ci{k-m-m}, 1.2.3)}
Let $X$ be a nonsingular variety and $\pi: X \to S$ be a proper
 morphism
onto a variety $S$. Assume that $N$ is a Cartier divisor on $X$
and that $M$ and
$\Delta$ are  $\rat$-Cartier divisors on $X$  with the following properties:
$(1)$  $M$ is $\pi$-nef and $\pi$-big,
$(2)$ the support of $\Delta$ is a divisor with normal crossings,
and $\lfloor \Delta \rfloor = 0$, 
and
$(3)$  $N\equiv M+ \Delta$.

\smallskip
\noindent
Then $R^i\pi_*\odix{X} (K_X + N) =0$ for $i>0$. 
\end{tm}

\begin{tm}
\label{inj}
{\rm (Cf. \ci{koshaf}, 10.13 and 9.17,  and \ci{e-v}, 5.12.b)}
Let $\pi:X \to S$ and $\Delta$ be as above with $X$ projective,
$D$ be an effective Cartier divisor on $X$ such that it does not dominate
$S$ via $\pi$,  and $L$ be a nef and big
$\rat $-Cartier divisor   on $S$. 
Let $N$ be a Cartier divisor such that 
 $N\equiv \Delta + \pi^*L$.

\noindent
Then the following  natural homomorphisms are injective for every $i\geq 0$:
$$
H^i (X, K_X +  N) \longrightarrow H^i(X,  K_X + N + D).
$$
\end{tm}

\medskip
\noindent
{\bf Local Seshadri constants.}
Good references for what follows are \ci{dem94}, \S7 and \ci{ekl}.
\begin{defi}
\label{seshco}
{\rm
Let $X$ be a complete variety, $L$ be a {\em nef} $\,$
$\rat$-Cartier divisor  on $X$, and $x$ be a point
on $X$. The following nonnegative real number
is called the {\em Shesadri constant of $L$ at $x$}:
$$
\e(L,x)= \inf \{ \frac{L\cdot C}{mult_x C} \},
$$
where the infimum is taken
over all integral curves passing through $x$ and $mult_x C$ is the 
multiplicity
of $C$ at $x$. 
} 
\end{defi}
 Let $x$ be a  point in $X_{reg}$,  $b_x: X' \to X$ be the blowing-up
of $X$ at $x$ and  $E$ be the corresponding exceptional divisor
on $X'$.
The $\rat$-Cartier divisor $b_x^* L$ on $X'$ is nef as well.
In particular,
there is a well-defined nonnegative real number:
$$
\e'(L,x) : = \sup \{\e' \in \rat | \, b_x^*L - \e' E \,\,\, is \,\, nef \,\}.
$$
It is clear that the 
$\real$-Cartier divisor 
$ b_x^*L - \e'(L,x) E$ is nef, 
and that the 
$\rat$-Cartier divisor
$ b_x^*L - \e' E$ is nef for every rational number $\e'$ with the property
that 
$0\leq \e' \leq \e' (L,x)$.

\noindent
\begin{fact}
\label{same}
 We have that $\e(L,x)=\e'(L,x)$ for every   $x\in X_{reg}$.
{\rm
This follows  from 
 the formula:
$(b_x^*L - \e E)\cdot \tilde{C}= L\cdot C - \e \, mult_x C$, where
$\e$ is any real number,  $\tilde{C}$ is any 
integral curve
in $X'$ not contained in $E$ and $C:=b_x(\tilde{C})$.
}
\end{fact}

\smallskip
We collect the simple properties of $\e(L,x)$ which, together with Theorem
\ref{ekltm}, are
essential in the sequel of the paper.

\begin{lm}
\label{basic} Let $L$ be a nef $\,\rat$-Cartier divisor on
a complete variety $X$ and $x$ be a  point in $X_{reg}$. Then 

\noindent
{\rm (\ref{basic}.1)}  $L^n \geq \e(L,x)^n$;

\noindent
{\rm (\ref{basic}.2)}  for every $t\in \rat^+$, $\e(tL,x)=t\,\e(L,x)$;
 
\noindent
{\rm (\ref{basic}.3)}  Let $f: X' \to X$ be a  proper and birational
morphism and $x$ be a point on $X$ over which $f$ is
an isomorphism; then $\e(L,x)=\e(f^*L, f^{-1}\{x\})$;

\noindent
{\rm (\ref{basic}.4)}  
if $L$ is Cartier, ample and generated by its global sections on $X$, then
$\e(L,x) \geq 1$;

\noindent
{\rm (\ref{basic}.5)} 
if $L$ is  Cartier and the global sections of $L$ generate jets of order 
$s$ 
at $x$,
{\em i.e.} the natural evaluation map
$H^0(X,L) \to \odixl{X}{L}/{\frak m}_x^{s+1} \odixl{X}{L}$ is surjective, 
then
$\e(L,x) \geq s$.

\end{lm}

\noindent
{\em Proof.} The first property follows  from the fact that
since $b_x^*L - \e(L,x)E$ is nef, then its top self-intersection is 
nonnegative.
The second one is an obvious consequence of the bilinearity of the 
intersection
product.

\noindent
The third property follows from the fact that there is a natural bijection,
given by taking strict transforms, 
between  the sets of integral curves on $X$ through $x$ and on $X'$ through
$x':=f^{-1}\{x\}$.
 If $C$ and $C'$ correspond to each other in this bijection, then
$L \cdot C= b_x^*L \cdot C'$ and $mult_xC=mult_{x'}C'$
so that the two local Seshadri constants are the same.

\noindent
If $L$ is ample, Cartier and generated by its global sections on $X$, 
then 
  the rational map $\varphi$  defined by
$|L|$ is a finite morphism. Let $C$ be any integral curve  on $X$ passing
through $x$. Since $C$ is not contracted by $\varphi$, there is an 
effective 
divisor $D$ in $|L|$ passing through $x$ but not containing $C$.
It follows that $L\cdot C=D\cdot C \geq mult_x C$. This implies
the third property.  

\noindent
Finally, if $s=0$, then  there is nothing  left to prove. Assume that 
$s\geq 1$.
Then the  global sections of $L\otimes {\frak m}_x^{s}$ 
generate $L\otimes {\frak m}_x^{s}$ at $x$. Given
any integral curve $\tilde{C}$ on $X'$ not contained in $E$, we find
a divisor $D \in |L\otimes {\frak m}_x^s|$ not containing the curve
$b_x (\tilde{C})$. It follows that the effective divisor $b_x^* (D) \in
|b_x^*L -sE|$ does not contain $\tilde{C}$. In particular,
$(b_x^*L -sE)\cdot \tilde{C} \geq 0$. This is enough to establish the 
last property.
\blacksquare

\medskip
If $X$ is complete
and $L$ is a nef $\rat$-Cartier divisor on $X$, then Shesadri's
criterion of ampleness asserts that $L$ is ample  iff
$\e(L): = \inf \{ \e(L,x) | \, x \in X \} >0$. 
An example of R. Miranda's
(cf. \ci{dem94}, 7.14) shows that given any positive real number
$\e$, there exists a nonsingular rational surface $X$, a point $x \in X$
and an ample 
line bundle $L$ on $X$, such that $\e(L,x) \leq \e$. 

\noindent
In particular,
we {\em cannot} expect to have  a statement of the 
form: {\em let $X$ be a nonsingular projective variety of dimension $n$ 
and $L$ be
an ample line bundle on $X$, then $\e(L) \geq C_n$, for some positive 
constant 
depending only on $n$}.

\smallskip
What is known is the following result of Ein, K\"uchle and Lazarsfeld. 
The authors prove it for  projective varieties, but 
by Chow's Lemma and  Lemma \ref{basic}.3  the statement is true for every 
complete variety.

\begin{tm}
\label{ekltm}
{\rm (Cf. \ci{ekl}, Theorem 1)}
Let $L$ be a nef and big Cartier divisor on a complete  variety $X$ 
of 
dimension $n$.  Then at a very general point $x$ on $X$ we have:
$$
\e (L,x) \geq \frac{1}{n}.
$$
\end{tm}

The example that follows   shows that
 Theorem \ref{ekltm} cannot hold as stated for an ample and effective 
 integral $\rat$-Cartier $\rat$-divisor
 on  a normal projective  variety.  
As is pointed out in \ci{ekl},
if $m$ is the smallest positive integer such that $mL$ is Cartier, then
Theorem \ref{ekltm} holds if we replace ``$\e(L,x) \geq \frac{1}{n}$" by
``$\e(L,x) \geq \frac{1}{nm}$."

\begin{ex}
\label{exindex}
{\rm
Let $S_m\subset \pn{m+1}$ be the surface which is a cone of vertex $v$ 
over the rational normal curve of degree $m$
in $\pn{m}$ and $\frak l$ be any line belonging to the ruling of $S_m$. 
The   Weil divisor
$\frak l$  is an integral  $\rat$-Cartier $\rat$-divisor; to be precise 
it is
 $m$-Cartier. The Cartier divisor $m{\frak l}$ is very ample
so that, for every $x \in S_m \setminus \{v \}$, we have that
$\e({\frak l},x) \geq \frac{1}{m}$. On the other hand,
fix $x \in S_m \setminus \{v \}$ and  let $ C$ be the line
on $S_m$ passing through $x$. We have ${\frak l} \cdot C=\frac{1}{m}$, so 
that
$\e( {\frak l},x ) \leq \frac{1}{m}$. It follows that 
$\e({\frak l},x) = \frac{1}{m}$, for every $x \in S_m \setminus \{v \}$.
}
\end{ex}

\section{An effective nonvanishing theorem}
\label{nonvan}
In this section we prove a nonvanishing theorem very similar to 
\ci{koebpf}, Theorem 3.2. 
While the statement   is clearly inspired by 
\ci{koebpf}, Theorem 3.2,  its simpler proof is inspired by
\ci{dps}, Lemma 3.21.

 The basic nonvanishing
and global generation at a generic point follow easily
from \ci{koebpf}, Theorem 3.2 (and in fact are slightly
weaker than this latter result). 
The ``multiple-points higher-jets" statements  do not follow directly
from the results in the literature.

\noindent
Let us   point out that
 Koll\'ar's  result implies  a  version
of Theorem \ref{effnonvan} with $x$  general,   p=1  and  s=0.
 However, one can use
this result in place of Koll\'ar's in proving the  effective 
base-point-freeness result  1.1 of \ci{koebpf}.

The advantages of Theorem \ref{effnonvan} are at least two. 

\noindent
The former
is  the simplicity of its proof
which consists of basic yoga and one blowing-up procedure. However,
we should stress
that more often than not this becomes an effective result if used 
in conjunction with
the non-trivial result 
Theorem \ref{ekltm}. 

\noindent
The latter is that it is a ``multiple-point higher-jets" effective result 
which, 
at least in principle, can be applied to prescribed points
and  can give more than  mere nonvanishing. For example, one can use
this result to obtain increased lower bounds  
of log-plurigenera (cf. \ci{koebpf}, \S4). We shall see some other 
applications
in the following sections.

\begin{rmk}
\label{logsmooth}
{\rm
Let $g:Y \to S$ be a proper morphism of varieties with $Y$ nonsingular and
$\Delta$ a divisor on $Y$ such that  $Supp(\Delta)$  has simple normal 
crossings.

\noindent
By virtue of generic smoothness,
there exists a largest Zariski-open dense subset $U$ of $S$ such that
(i) $g_{|g^{-1}(U)}: g^{-1}(U) \to U$ is smooth, 
(ii) for every point $x \in S$, any irreducible component $F$ of the fiber
$F_x$ of $g$  is not contained
in $Supp(\Delta)$
and (iii) $Supp(\Delta)$ has simple normal crossings on $F$.
}
\end{rmk}

\begin{tm}
\label{effnonvan}
{\rm ({\bf Effective Nonvanishing})}
Let the following data be given.

\noindent
{\rm (\ref{effnonvan}.1)} 
  $(Y,\Delta)$:  a log-pair, where  $Y$  is nonsingular and complete,
$\lfloor\Delta \rfloor = 0$ and
 $Supp(\Delta)$ has simple normal crossings. 

\noindent
{\rm (\ref{effnonvan}.2)}
 $N$:  a Cartier divisor on
$Y$.

\noindent
{\rm (\ref{effnonvan}.3)} 

\noindent
-  $g: Y \to S$:   a proper  and 
surjective morphism 
onto    a complete  variety $S$ of positive dimension, 

\noindent
- $U=U(g,\Delta)$: 
the Zariski-open dense set of $S$ defined in {\rm  Remark \ref{logsmooth}},

\noindent
- 
$V$:  the Zariski-open dense
subset of $S$ over which the formation of
$g_*$ for $N$ commutes with base extensions.

 \noindent
{\rm (\ref{effnonvan}.4)}

\noindent
- 
$p$:  a positive integer, 

\noindent
- $\{s_1, \ldots, s_p\}$:   a $p$-tuple
of non-negative integers,

\noindent
-  $\{x_1, \ldots, x_p\}$:   $p$ distinct points in $U\cap V$.
 
\noindent
{\rm (\ref{effnonvan}.5)} 
$M$:  a $\rat$-Cartier divisor on $Y$ such that either 
it  is nef and  $g$-big, or $X$ is projective and $M\equiv 0$.

\noindent
{\rm (\ref{effnonvan}.6)} 
 $L_1, \ldots, L_p$: $p$  $\rat$-Cartier divisors   on $S$ such that
 all $L_j$ are nef and big
and either  

\noindent
(a) $\e(L_j, x_j) > \dim S + s_j$, $\forall j=1, \ldots, p$,
or 

\noindent
(b) $\e(L_j, x_j) \geq \dim S + s_j$, $\forall j=1, \ldots, p$ and 
$L_{j_0}^{\dim{S}} > \e(L_{j_0})^{\dim{S}}$ for at least one index
$j_0$, $1 \leq j_0 \leq p$.

\smallskip
\noindent
Assume that
$$
N \equiv K_Y + \Delta + M +  g^*\sum_{j=1}^p L_j.
$$

\smallskip
\noindent
Then the following natural map is surjective
$$
H^0(X,N) \simeq H^0(S, g_*N) \surj
\bigoplus_{j=1}^p  \, \, \frac{g_*N}{{\frak m}_{x_j}^{s_j+1} \cdot  g_*N}.
$$
In particular, if $g_*N$, which is torsion-free,
 is not the zero sheaf, then $H^0(X,N)\not= \{0\}$.
\end{tm}

\begin{rmk}
\label{why}
{\rm 
The reason for calling this theorem ``Effective Nonvanishing" is the 
last assertion of the theorem and  the fact that, for example, 
if all the $L_j$
were Cartier, then we  could make sure, 
by virtue of
Theorem \ref{ekltm}, that condition
(\ref{effnonvan}.6) is fulfilled  at very general points
by taking sufficiently high multiples of the $L_j$.

\noindent
Note also that the conclusion of the theorem holds trivially also for 
$\dim{S}=0$, but that
in this case (2.1.6) is not meaningful.
}
\end{rmk}

\noindent
{\em Proof.} 
The proof is divided into two cases. The former deals with $M$
nef and $g$-big. The latter with $X$ projective and $M\equiv 0$.
Each  case is divided into two sub-cases corresponding to the 
 two numerical assumptions (a) and (b) in (\ref{effnonvan}.6).

\smallskip
\noindent 
CASE I: {\em $M$ is nef and $g$-big}.

\noindent
First we show  that in this case $U=U\cap V$.
By virtue of \ref{van}, we know that
$R^ig_* N=0$ for $i >0$.  By the smoothness  of $g$ over $U$, 
$N$ is flat over $U$. By well-known results of Grothendieck
(see \ci{gro}, III.7.7.10) 
$g_*N$ is locally free on $U$ and  the formation
of $g_*$ commutes with base extension over $U$. 

\smallskip
\noindent
In particular,
if $Y_{x_j}^{\s}:= Y\times_S  Spec \, (\odix{S,x_j}/{\frak m}_{x_j}^{\s})$
is the ``$\s$-thickened fiber" of $g$ at $x_j$ and
$N_{x_j}^{\s}$ is the pull-back 
 of $N$ to $Y_{x_j}^{\s}$ via the natural projection, then the following
 natural maps are  isomorphisms:
$$
\frac{g_*N}{{\frak m}_{x_j}^{s_j+1} \cdot  g_*N}=
g_*N \otimes ( \odix{S,x_j}/{\frak m}_{x_j}^{s_j+1} ) \longrightarrow
 H^0(Y_{x_j}^{s_j+1},
 N_{x_j}^{s_j+1}).
$$

\smallskip
\noindent
To  prove CASE I  it is enough to show that
the natural map
\begin{equation}
\label{1}
H^0(Y,N) \longrightarrow 
\bigoplus_{j=1}^p \, \, H^0(Y_{x_j}^{s_j+1}, N_{x_j}^{s_j+1}),
\end{equation}
which factors through $g_*N \otimes \odix{S,x_j}/{\frak m}_{x_j}^{s_j+1}$,
is  surjective.

\smallskip
\noindent
Consider the following cartesian diagram:
$$
\begin{array}{lll}
 \hspace{1cm} Y' &
 \stackrel{B}\longrightarrow & 
Y   \\
\hspace{1cm} \downarrow {g'} & \ & \downarrow g \\
  \hspace{1cm} S'            & \stackrel{b}\longrightarrow & 
S               
\end{array}
$$
where $b$ is the blowing-up of $S$ at all the simple points $x_j$. Let
$F:=\coprod F_j$ be the scheme-theoretic-fiber of $g$ corresponding to
the union of the points $x_j$, $j=1,\ldots, p$.
Since $g$ is smooth
over $U$ and all the $x_j$ are in $U$, we see that $B$ coincides with the 
blowing-up
of $Y$ along $F$. In particular, $Y'$ is a nonsingular variety.
 Let $E=\sum E_j$ be the exceptional divisor of $b$ and $D=\sum D_j$
the one of $B$; we have that 
$D_j={g'}^*E_j$, for every $j=1, \ldots, p$. 

\smallskip
\noindent
The map (\ref{1}) is surjective iff
the natural map
$H^1(Y',B^*N-\sum (s_j+1)D_j) \to  
H^1(Y',B^*N)$ is injective.  
It is this injectivity that we are 
going to establish  using 
Theorem \ref{van}.  

\smallskip
\noindent
Note that $K_{Y'} \approx B^*K_Y + (\dim S -1)\sum{D_j}$ and
 that 
since no irreducible component of any $F_j$ is contained 
in any $\Delta_i$ and if any such component meets any
$\Delta_i$ it does so transversally, we have that 
a) $\Delta':=B^*\Delta=B^{-1}_*\Delta$, {\em i.e.}
the pull-back is the strict transform, 
 b) $\lfloor \Delta' \rfloor=0$ and  c) the support
of $\Delta'$ has simple  normal crossings. The following numerical equality
 is easily checked:
\begin{equation}
\label{2}
B^*N - \sum{(s_j+1)D_j} \equiv K_{Y'} + \Delta' + B^*M + B^*g^*\sum{L_j} -
\sum{(\dim{S} + s_j)D_j}. 
\end{equation}

\smallskip
\noindent
SUB-CASE I.A:  {\em Assume that $\e(L_j,x_j)> \dim{S} + s_j$, for 
every index $j$, $1 \leq j \leq p$.}

\noindent
Since  for every index $j$ we have that  $\e(L_j,x_j)> \dim{S} + s_j$, 
there exists
a positive rational number $0<t<1$ such that
$\e((1-t)L_j,x_j)> \dim{S} + s_j$ for every $j$, 
$1\leq j \leq p$. Using the fact that $B^*g^*={g'}^*b^*$ 
we can re-write
the r.h.s. of equation (\ref{2}) as
\begin{equation}
\label{3}
K_{Y'} + \Delta' + B^*(M + tg^*\sum{L_j}) + {g'}^*\sum
{\left[ b^*(1-t)L_j - 
(\dim{S}+ s_j)E_j \right]}.
\end{equation}
The last summand is nef by the very definition of $\e((1-t)L_j,x_j)$.

\noindent
Since $M$ is nef and $g$-big and $t>0$, the $\rat$-divisor $M + 
tg^*\sum{L_j}$ 
is nef and big. In particular, $B^*(M + tg^*\sum{L_j})$ is nef and big.
It follows that the  l.h.s. of (\ref{2}) is a Cartier divisor satisfying
the assumptions of Kawamata-Viehweg Vanishing Theorem  so that 
$H^1(Y', B^*N-\sum (s_j+1)D_j)=\{0\}$ and (\ref{1}) is surjective. 

\medskip
\noindent
SUB-CASE I.B: {\em Assume that $\e(L_j, x_j) \geq \dim S + s_j$, $\forall 
j, \,
1\leq j \leq  p$ and that  
$L_{j_0}^{\dim{S}} > \e(L_{j_0})^{\dim{S}}$ for at least one index
$j_0$, $1 \leq j_0 \leq p$.}

\noindent 
Using the fact that $B^*g^*={g'}^*b^*$  and isolating the index $j_0$ we 
write
the r.h.s. of (\ref{2}) as
\begin{equation}
\label{4}
K_{Y'} + \Delta' + B^*M + \sum_{j\not= j_0}{{g'}}^*
\left[
b^*{L_j}- (\dim{S}+ s_j)E_j\right] +
{{g'}}^*
\left[ b^*{L_{j_0}}- (\dim{S}+ s_{j_0})E_{j_0} \right].
\end{equation}
Since $M$ is nef and $g$-big and $\sum_{j\not= j_0}{{g'}}^*
(b^*{L_j}- (\dim{S}+ s_j)E_j)$ is nef we see that 
$B^*M + \sum_{j\not= j_0}{{g'}}^*
(b^*{L_j}- (\dim{S}+ s_j)E_j)$
 is nef and ${g'}$-big. 
Since $L_{j_0}^{\dim{S}} >  \e(L_{j_0},x_{j_0})^{\dim{S}}$, we see, as in
the proof of 
Lemma \ref{basic}.1, that 
$
(b^*{L_{j_0}}- (\dim{S}+ s_{j_0})E_{j_0})$
is nef and big. It follows that 
$ B^*M + \sum_{j\not= j_0}{{g'}}^*
(b^*{L_j}- (\dim{S}+ s_j)E_j) +
{{g'}}^*
(b^*{L_{j_0}}- (\dim{S}+ s_{j_0})E_{j_0})$
is nef and big  and we conclude as in SUB-CASE I.A.

\medskip
\noindent
CASE II: {\em $X$ is projective, $M\equiv0$ and the points $x_j$ are in 
$U\cap V$.}

\noindent
We by-pass the first paragraph in the proof of CASE  I. We proceed {\em 
verbatim}
as in that case until we hit again (\ref{2}).
We delete $M$.
We can again divide the analysis into two separate sub-cases. We do so and
obtain that in the two distinct sub-cases
the l.h.s. of (\ref{2}) 
is  numerically equivalent to the r.h.s. of (\ref{3})
and (\ref{4}), respectively
and, in both cases,  we are in the position to apply Theorem \ref{inj} to
the morphism ${g'}:Y' \to S'$ and infer the desired
injectivity statement.
\blacksquare

\section{Applications}
\label{firstapp}
The local Seshadri constant
 can be linked, via Kawamata-Viehweg Vanishing Theorem to the production
of sections for the adjoint to nef and big line bundles.   
This observation is due to Demailly; see \ci{dem94}, Proposition 7.10
 and \ci{ekl}, Proposition 1.3.  
In this section we apply Theorem \ref{effnonvan} to nef vector bundles.
Actually,  a factor $\det{E}$ appears and is necessary in our proof. We 
ignore if it 
is necessary for the truth of the  various statements  that follow.
First we fix some notation.

\bigskip
Let $E$ be a rank $r$ vector bundle on a nonsingular
complete variety  $X$. We denote by ${\Bbb P}_X(E)$ the projectivized bundle
of hyperplanes,
by $\pi: {\Bbb P}_X(E) \to X$ the natural morphism and  by $\xi$ or $\xi_E$
the tautological line bundle $\odixl{{\Bbb P}_X(E)}{1}$. We say that
$E$ is nef if $\xi$ is nef.

Let $p$ be any
 positive integer.
We say that {\em the global sections of $E$ generate 
jets of order
$s_1,\ldots, s_p \in {\Bbb N}$ at  $p$ distinct  points
 $\{ x_1, \ldots ,x_p \} $ of $X$}
if the following natural map is  surjective:
$$
H^0(X,E) \longrightarrow \bigoplus_{\i=1}^p  E_{x_{i}} 
\otimes \odix{X}/{\frak m}^{s_{i}+1}_{x_{i}}.
$$

\noindent
We say that 
{\em 
the global  sections of $E$ separate 
$p$  
distinct points $\{x_1, \ldots, x_p \}$
 of  $X$
} 
if the above holds with
all $s_{i}=0$. The case $p=1$ is equivalent to $E$ being generated by 
its global sections ({\em generated}, for short) at the point in question.

\smallskip
\noindent
{\bf Rational maps to Grassmannians.}
 Let $V:= H^0(X,E)$
and $h^0:=$ $h^0(X,E):=$ $\dim_k H^0(X,E)$.
Consider  the Grassmannian $G:= G(r,h^0)$ of $r$-dimensional quotients
of $V$,  the universal quotient bundle $\QQ$ of $G$ 
and  the determinant of $\QQ$, $\q$.

\noindent
As soon as $E$ is generated  at some point of $X$, 
we get a rational map
$\varphi:X --> G$ assigning to each point $y \in X$ where $E$ is 
generated 
the quotient $E_y \otimes k(y)$. 

\noindent
If $E$ is generated 
at every point of $X$, then $f:=\varphi$ is a morphism
and   $E \simeq f^* \QQ$. 

\noindent
It is clear that:

- $V$ separates arbitrary pairs of points of $X$ iff 
 $f$ is  bijective birational onto its image;

- If $V$ separates every pair of points of $X$ and generates
jets of order $1$
at every point of $X$, then  $f$ is a closed embedding (the converse 
maybe false 
if $r>1$).

\medskip
In the three propositions that follow we   generalize to the case of 
higher rank
 results in  \ci{ekl}. 
The analogues to these facts involving arbitrary $p$ and $\{s_1, \ldots, 
s_p\}$ 
are clear, and left to the reader. 
 We give the reference to the analogous 
results for line bundles, but  we prove only the first of the  three 
propositions to illustrate the method.

\noindent
\begin{pr}
\label{sect-s}
{\rm (Cf. \ci{ekl}, 1.3 and  4.4)}
Let $X$ be a nonsingular complete variety of dimension $n$. Let
   $E$  be a rank $r$ nef vector bundle on $X$,  $L$ be 
a nef and big $\rat$-Cartier divisor on $X$,  $\Delta'$ be a 
$\rat$-Cartier divisor
on $X$
such that $\lfloor \Delta' \rfloor =0$ and $Supp(\Delta')$ 
has simple normal crossings,
and $N'$ be a Cartier divisor on $X$ such that $N'\equiv L + \Delta'$.

\noindent
Let  $s$ be   a nonnegative integer and $x$ be a point of $X\setminus 
Supp(\Delta)$.

\noindent
Assume that 
either
$\e(L,x) > n+s$, or $\e(L,x)\geq n+s$ and $L^n > \e(L,x)^n$.

\noindent
Then $H^0(X, K_X \otimes E \otimes \det{E} \otimes N')$ generates 
$s$-jets at $x$
and the rational map $\varphi$ as above is defined.
Moreover,
$$
h^0(X, K_X \otimes E \otimes \det{E} \otimes N') \geq
r B(n+s,s). 
$$
In particular, if
$\cal L$ is  a  nef and big
Cartier divisor on $X$, then
$$
H^0(X, K_X \otimes E \otimes \det{E} \otimes {\cal L}^{\otimes m}) \geq
r B(n+s,s), \quad \forall \, m \geq n^2 +ns.
$$
\end{pr}

\noindent
{\em Proof.} 
Set 
$Y:={\Bbb P}_X(E)$, $S:=X$, $g:=\pi$, $\Delta:=g^*\Delta'$,
$M:=(r+1)\xi$, $N:=K_Y + (r+1)\xi + g^*N'$, $p=1$, $s_1=s$.
Note that  $M$ is nef and $g$-big and  that 
$g_*N=K_X \otimes E \otimes \det{E} \otimes N'$.

\noindent
Apply Theorem \ref{effnonvan}.
The only issue is whether $x\in U$; this is why the point $x$ is assumed 
to be 
outside of $Supp(\Delta)$.

\noindent
The lower bound on $h^0$ stems from the  surjection given by Theorem
\ref{effnonvan} and the fact that
$$\dim_k \odix{X,x}/{\frak m}_x^{s+1}=B(n+s,n).$$

\noindent
The statement about $\cal L$ is a special case after Theorem \ref{ekltm}:
there exists $x\in X$ such that $\e({\cal L},x)\geq 1/n$. If $m\geq n^2 
+ns$, then
$\e(m{\cal L}, x) \geq n+s$ and equality holds 
iff $\e({\cal L},x)=1/n$ and $m=n^2+ns$;
in this  case the inequality ${\cal L}^n \geq 1 > \e({\cal L},x)^n$
is automatic.
\blacksquare

\begin{pr}
\label{genbirat}
{\rm (Cf. \ci{ekl}, 4.5)}
Let $X$, $n$, $E$, $L$, $\Delta'$ and $N'$ be as above. 
Assume that either $n\geq 2$ and  $\e(L,x)\geq 2n$ for every
$x$ very general, or that $n=1$ and $\deg{N'}\geq 3$.

\noindent
Then 
the rational map  $\varphi$ associated with 
$H^0(X, K_X \otimes E \otimes \det{E} \otimes N' )$ is defined
and is birational onto its image.

\noindent
In particular, if ${\cal L}$  is a nef and big  Cartier divisor on $X$, then
the rational map  $\varphi$ associated with 
$H^0(X, K_X \otimes E \otimes \det{E} \otimes {\cal L}^{\otimes m} )$ 
is defined and birational onto its image for every $m \geq 2n^2$.
\end{pr}

\begin{pr}
\label{lt}
{\rm (Cf. \ci{ekl}, 4.6)}
Let $X$ be a complete variety  of dimension $n$ 
with only terminal singularities  and of global index $i$ such that
$K_X$ is nef and big, ({\em i.e.} $X$ is normal, $\rat$-Gorenstein and a 
minimal
variety of general type, and $i$ is the smallest positive integer such that
the Weil divisor class $iK_X$ is a Cartier divisor class), and
$E$ be a nef vector bundle
on $X$. 

\noindent
Then
the rational map associated with
$H^0(X, {\cal O}_X(miK_X) \otimes E \otimes \det{E}  )$ is 
defined and is birational onto its image for 
every $m\geq 2n^2 +1$.
\end{pr}

\medskip
The following follows
from results in \ci{koshafinv}, \S8. As is already 
pointed out in \ci{ekl},  a generically large algebraic fundamental group
on the base variety $S$ can be used to produce section by increasing the 
local Seshadri constants on  finite \'etale covers of $S$.
 The reader can consult
\ci{koshafinv} for the relevant definitions. 

\begin{pr}
\label{nonvanfundgrp}
{\rm (Cf. \ci{koshafinv}, 8.4)}
Let $X$ be a  normal and complete variety,  
$N'$ be an integral   big $\rat$-Cartier $\rat$-divisor on $X$, 
 and $E$ be a nef 
vector bundle
 on $X$. 

\noindent
Assume that $X$ has generically large algebraic fundamental group.

\noindent
Then $h^0(X, \odixl{X}{K_X +N'} \otimes E \otimes \det{E}  ) > 0$.
\end{pr}

\noindent
{\em Sketch of proof.} By the proof of \ci{koshafinv}, Corollary  8.4
and by the first part of the proof of \ci{koshafinv}, Theorem 8.3 we are 
reduced to the case in which $X$ is nonsingular and $N'\equiv L + 
\Delta$, where
$L$ and $\Delta'$ are   $\rat$-Cartier divisors, 
$L$ is  nef and big,  $\lfloor \Delta' \rfloor =0$ and $Supp(\Delta')$ 
has 
simple normal crossings.

\noindent
Pick a point $x\in X$ such that $\e(L,x) >0$. By \ci{koshafinv}, Lemma 8.2
there is a finite \'etale map of varieties $m:X'' \to X$ and a point 
$x''\in X''$
such that $\e(m^*L,x'') \geq n$. 

\noindent
Denote $\deg{m}$ by $d$,  $m^*L$ by $L''$, $m^*\Delta'$ by $\Delta''$,
 $m^*N'$ by $N''$
 and $m^*E$ by $E''$.

\noindent
Apply Proposition \ref{sect-s} to $X''$, $L''$, $\Delta''$, $N''$,  $E''$ 
and $s=0$.
We get $h^0(K_{X''}\otimes E'' \otimes \det{E''} \otimes N'')>0$.

\noindent
Kawamata-Vieheweg Vanishing Theorem applied to the nef and big $\rat$-divisor
$(r+1)\xi_{E''} + {\pi''}^*L''$ gives, via Leray spectral sequence,
$h^i(X'', K_{X''}\otimes E'' \otimes \det{E''} \otimes N'')=0$, for every 
$i>0$.
The analogous statement holds  on $X$.

\noindent
The above vanishing and the multiplicative behavior of  Euler-Poincar\'e
characteristics of coherent sheaves under finite \'etale maps of nonsingular
proper varieties gives:
$$
h^0(X, K_{X}\otimes E \otimes \det{E} \otimes  N')=
\chi (X,-)=
\frac{1}{d}\chi (X'',-'')=\frac{1}{d}
h^0(X'',K_{X''}\otimes E'' \otimes \det{E''} \otimes N'') >0.
$$
\blacksquare

Let us point out a consequence of $\ci{koshafinv}, 8.3$ as a corollary 
to the result above.
Recall that the integers  
$I_{lm}^i:
=h^i(X,S^l (\Omega_X^1) \otimes K_X^{\otimes m})$
are birational invariants of a  nonsingular and complete  variety $X$ 
for every $m,l\geq 0$ and that they are independent
of the more standard invariants like the plurigenera
or     the cohomology groups of the sheaves 
$\Omega^p_X$; for some facts about these invariants and some references
see \ci{rack}.  The assumptions 
of the ``sample"  corollary  that follows are fulfilled, for example, by 
projective varieties  whose universal covering space is the unit ball
in $\comp^n$.

\begin{cor}
\label{eximple}
{\rm (Cf. \ci{koshafinv}, 8.5)}
Let $X$ be a nonsingular complete variety with $K_X$ nef and big, 
$\Omega^1_X$
nef and generically large algebraic fundamental group.

\noindent
Then $I_{1m}^0\geq 0$ for every $l \geq 0$ and $m\geq 3$. 
\end{cor}

\medskip
We now observe that the global generation results of
Anghern-Siu, Demailly, Tsuji and Siu can be used to
deduce analogous statements for vector bundles of the form
$K_X^{\otimes a} \otimes E \otimes \det{E} \otimes L^{\otimes m}$, where
$E$ and $L$ are   a nef vector bundle and an ample line bundle on $X$, 
respectively. The idea is simple: 
once the  sections of a
line bundle of the form  ${\cal L}:=K_X+mL$ 
generate the $s$ jets  at every point,
 the  local Seshadri constant is at  least $s$ at every point by virtue 
of Lemma
\ref{basic}.5. We then use Proposition \ref{sect-s}.
However, this idea is applied here  in a rather primitive
way; we expect these results to be far from otpimal.

\noindent
We shall give  statements  concerning
$p=1,2$ and low values for the jets. In the same way one can prove statements
concerning more points and higher jets. 
We omit the details.

For ease of reference we collect the line bundle results in the literature
in the following result. First some additional notation.
Let $n$  and $p$ be  positive integers and $\{ s_1, \ldots , s_p \}$
be a $p$-tuple of nonnegative integers. Let us
define the following integers:
$$
m_1 (n,p) : = \frac{1}{2}(n^2 +2pn -n +2 ),
$$
$$
m_2(n,p;s_1, \ldots , s_p)= 2n \sum_{i=1}^p B(3n + 2s_i -3, n) + 2pn +1.
$$

\begin{tm}
\label{effres}
Let $X$ be   a nonsingular  projective variety of dimension
$n$,  and $L$  be an  ample 
Cartier divisor on $X$.

\medskip
\noindent
{ \rm (\ref{effres}.1) (Cf. \ci{siu94b})} if $m\geq m_2(n,p; s_1, \ldots, 
s_p)$, 
then
the global sections of $2K_X + mL$
generate simultaneous jets of order
$s_1,\ldots, s_p \in {\Bbb N}$ at arbitrary  $p$ distinct points of $X$;

\medskip
\noindent
{ \rm (\ref{effres}.2) (Cf. \ci{an-siu})} If $m\geq m_1(n,p)$, then
the global sections of $K_X + mL$ separate arbitrary
$p$ distinct points of $X$;

\end{tm}

\begin{tm}
\label{effresvb}
Let $X$, $n$ and $L$ be as above. Let $E$ be a nef vector bundle on $X$.
Then the vector bundles $K_X^{\otimes a} \otimes E \otimes \det{E} \otimes
L^{\otimes m}$:

\noindent
{\rm (\ref{effresvb}.1)}
 are generated by their global sections and the associated morphism 
to a Grassmannian $f:X \to G$ is finite, for $a=2$ and for every
$m\geq (1/2)(m_2(n,1;2n)+1)$;

\noindent
{\rm (\ref{effresvb}.2)} 
have global sections which separate arbitrary pairs of points,
$1$-jets at an arbitrary point, and $f$ is a closed embedding, for
$a=2$ and for every $m\geq (1/2) ( m_2(n,1;4n) + 1)$;

\noindent
{\rm (\ref{effresvb}.3)}
 are generated by their global sections and 
 $f$ is finite, for $a=n+1$ and for every
$m\geq nm_1(n,1)$;

\noindent
{\rm (\ref{effresvb}.4)} 
have global sections which separate arbitrary pairs of points,
 $1$-jets at an arbitrary point, and $f$ is a closed embedding, for
$a=2n+1$ and for every $m\geq  2nm_1(n,1)$.

\end{tm}

\noindent
{\em Proof.}
Let us observe that all the vector bundles in question are ample.
One sees this easily by observing that $K_X + (n+1)L$
is always nef (Fujita) and that ``nef $\otimes$ ample $=$ ample." 
 As soon as 
$f$ is defined, these bundles are pull-backs under $f$ so that they 
can be ample only if $f$ is finite.

\smallskip
\noindent
Let $L':= K_X + (1/2) (m_2(n,1;2n) + 1)L$. By  virtue of
Theorem \ref{effres}.1, the global sections of
$2L'$  generate $2n$ jets at every point $x\in X$. By virtue of
Lemma \ref{basic}.5, 
$\e(L',x)\geq n$ for every $x\in X$. We can apply Proposition \ref{sect-s}
which assumptions are readily verified.  This proves
(\ref{effresvb}.1).

\noindent
The proof of (\ref{effresvb}.2) is similar. 
We observe that we need $\e(L',x) \geq 2n$
to separate points and $\e(L',x)\geq n+1$ to separate $1$-jets. We then use
Proposition \ref{genbirat} in the former case and Proposition \ref{sect-s}
with $s=1$
in the latter.

\noindent
 (\ref{effresvb}.3) and (\ref{effresvb}.4) are proved similarly using
Theorem  \ref{effres}.2 and  Lemma \ref{basic}.4.
\blacksquare

\section{Better bounds for global generation}
\label{einz}
In this section we greatly improve upon 
Theorem \ref{effresvb}.1 and \ref{effresvb}.3.
The method is similar to the one of the previous section. However, it does not 
use Seshadri constants. It needs a similar,  local positivity result which allows
one to apply the same techniques used before, and based on 
Theorem \ref{effnonvan}, to produce sections. Once the local positivity at one point
has been established,   the technique employed
in \ref{effnonvan} emerges in all its simplicity.

Let us recall, for the readers's convenience few basic facts
 about  the algebraic counterparts to Nadel Ideals.
The reference is \ci{ein}.

\smallskip
Let $X$ be a nonsingular projective variety and $D$ be an effective $\rat$-divisor.
Let $f: X'\to X$ be an embedded resolution for $(X,D)$. The integral divisor
$K_{X'/X} - f^* \lfloor D \rfloor$ can be written as $P-N$, where $P$ and $N$ are integral divisors without common components and $P$ is $f$-exceptional.

The {\em multiplier ideal} ${\cal I} (D)$
 associated with $(X,D)$ is, by definition, 
$$
 {\cal I} (D) := f_* \odixl{X'}{P-N} = f_* \odixl{X'}{-N}
\subseteq \odix{X}.
$$
One checks that this ideal sheaf is independent of the resolution chosen and that $\odixl{X'}{P-N}$ has trivial higher direct images.
As a consequence, we get the following vanishing result.

\begin{pr}
\label{evr}
{\rm (Cf. \ci{ein}, 1.4)} Let $X$ be a nonsingular projective variety, $L$ be a line bundle on $X$ and $D$ be an effective $\rat$-divisor on $X$.
Assume that $L-D$ is nef and big.

\noindent
Then $H^j(X, K_X \otimes L \otimes {\cal I} (D)) = \{0 \}$, for
every $j>0$.
\end{pr}
The following lemma is an easy consequence of the definitions and is a
functorial property of these ideals.
\begin{lm}
\label{funct}
Let $\pi : P \to X$ be a smooth morphism of nonsingular projective varieties
and $D$ be an effective $\rat$-divisor on $X$.
Then $\pi^* {\cal I} (D) = {\cal I} (\pi^* D)$.
\end{lm}
{\em Proof.} Consider the following cartesian diagram
$$
\begin{array}{ccc}
P' 
&
\stackrel{f'}\longrightarrow       
& 
P
\\
\downarrow{\pi'}
&
{\Box}
&\downarrow{\pi}
 \\
X'& \stackrel{f}\longrightarrow         & X
\end{array}
$$
where $f:X' \to X$ is an embedded resolution of singularities of the log-pair
$(X,D)$. Since $\pi$ is smooth, $f': P' \to P$ is a resolution of $(P, \pi^*D)$.

\noindent
We have ${\cal I}(\pi^* D) = f'_* (K_{P'/P} - \lfloor  {f'}^* (\pi^* D)\rfloor  )=$
$ {f'}_* ({\pi '}^* K_{X'/X} - 
\lfloor   {\pi '}^* f^*D\rfloor  )=$ ${f'}_* ( {\pi '}^* K_{X'/X} - {\pi '}^* \lfloor  f^*D\rfloor  )=$
${f'}_* ( {\pi '}^* ( K_{X'/X} - \lfloor  f^*D\rfloor  )) =
$
$\pi^* ( f_* (K_{X'/X} - \lfloor  f^*D\rfloor  )) = $ $\pi^* {\cal I} (D)$,
where: the second equality holds because the formation of the sheaf of relative differentials $\Omega^1_{*/*}$ commutes with base change and the relative cononical sheaf is, in the simple case under scrutiny, the determinant of $\Omega^1_{*/*}$; the third equality holds because $\pi '$ is smooth; the fifth stems from the fact that cohomology commutes with the flat base extension
$\pi$.
\blacksquare

The following result is a
$\rat$-divisors reformulation of the result of Anghern-Siu.
The result is due to Koll\'ar \ci{koslp}. The formulation
given below in terms of algebraic multiplier ideals is due to
Ein \ci{ein}.

\begin{tm}
\label{ein}
Let $X$ be a nonsingular projective variety of dimension $n$ and 
$L$ be an ample line bundle on $X$ such that
$$
L^d \cdot Z > B(n+1,2)^d
$$
for every $d$-dimensional integral cycle $Z$ on $X$.
Then, for every point $x\in X$ there exists an effective $\rat$-divisor 
$D$ such that $D\equiv \lambda L$ for some positive rational number
$0<  \lambda <1$ and $x$ is in the support of an isolated
component of $V( {\cal I} (D) )$.
\end{tm}
\begin{rmk}
\label{sevpts}
{\rm
A similar statement holds if we consider several distinct points.
}
\end{rmk}

\begin{tm}
\label{asvb}
Let $\pi : P \to X$ be a smooth morphism with connected fibers
of nonsingular projective varieties,
$n$ be the dimension of $X$, $M$ be a nef and $\pi$-big line bundle
on $P$, $\cal L$ be an ample line bundle on $X$ such that
$$
{\cal L}^d \cdot Z > B(n+1,2)^d
$$
for every $d$-dimensional integral cycle $Z$ on $X$.

\noindent
Then, the vector bundle $\pi_* (K_P + M)\otimes {\cal L}$ is generated by its global sections.

\noindent
In particular, if $L$ is any ample line bundle on $X$, then we can choose,
${\cal L}:= B(n+1,2) \, L$.
\end{tm}
{\em Proof.}
Let $x \in X$ be an arbitrary point. 
Let $D$ be a $\rat$-divisor as
in Theorem \ref{ein}.
Since ${\cal L} - D$ is ample and $M$ is nef and $\pi$-big, the $\rat$-divisor
$M+ \pi^* ({\cal L} -D)$ is nef and big on $P$.
The smoothness of $\pi$ implies, by virtue of Lemma \ref{funct}, that
$\pi^* {\cal I} (D) = {\cal I} (\pi^* D)$. It follows that
$H^1(P, (K_P + M  + \pi^* {\cal L}) \otimes \pi^* {\cal I} (D))=$
$ H^1(P, (K_P + M + \pi^*{\cal L}) \otimes  {\cal I} (\pi^* D) ) =\{0\}$,
 the second equality stemming from
Ein's version of Nadel Vanishing Theorem Proposition \ref{evr}.

Since $V({\cal I} (D))$ has isolated support at $x$, we conclude that, if we denote by $F_x$ the fiber of $\pi$ over $x$:
$$
H^0(P, K_P + M + \pi^* {\cal L} ) \surj 
H^0 (F_x,   (K_P + M + \pi^* {\cal L})\otimes \odix{F_x} ).
$$
The result follows by the natural identification between the map given above and
the map
$$
H^0 (X, \pi_*(K_P + M) \otimes {\cal L}) \to \pi_* (K_P+M) \otimes {\cal L}
\otimes \odix{X}/{\frak m}_x,
$$
which holds because $R^1\pi_* (K_P+M)=0$ is 
the zero sheaf by relative vanishing.
\blacksquare
\begin{cor}
\label{good}
Let $X$ be a nonsingular projective variety of dimension $n$, $E$ be a nef vector bundle on $X$, and $\cal L$ be an  ample line bundle
on $X$.
Assume that
$$
{\cal L}^d \cdot Z > B(n+1,2)^d
$$
for every $d$-dimensional integral cycle $Z$ on $X$.
\noindent
Then $K_X \otimes E \otimes \det{E} \otimes {\cal L}$ is generated by its global sections at every point of $X$.

\noindent
In particular, if $L$ is any ample line bundle on $X$, then we can choose
${\cal L} = B(n+1, 2) \, L$.
\end{cor}
{\em Proof.} 
Set $P:= {\Bbb P }(E)$, $\pi:=$ the natural projection onto $X$, 
$M: = (r+1)\xi_E$, where $r$ is the rank of $E$, and apply
Theorem
\ref{asvb}.
\blacksquare
\begin{rmk}
{\rm A similar statement hold for the simultaneous generation at several points;
see Remark \ref{sevpts}. The same is true for Theorem \ref{asvb}. 
}
\end{rmk}
\begin{rmk}
{\rm
The paper \ci{deshm} contains similar results without the factor
$\det{E}$. However, the assumption $E$ nef is there replaced by a stronger curvature condition on $E$, and the methods are purely analytic
}
\end{rmk}
\begin{rmk}
{\rm 
We do not know if similar statements hold without the factor $\det {E}$.
}
\end{rmk}

\bigskip
\noindent
 $1991$ {\em Mathematics Subject Classification}: 14E20, 14F05, 14J60, 
14Q20. 

\smallskip
\noindent
{\em Keywords and phrases}:  Effectivity, algebraic fundamental group,
local Seshadri constant, maps to Grassmannians,
 nef  vector bundle,
non-vanishing theorems, effective global generation, jets.

\bigskip
\noindent
Author's address:
Max-Planck-Institut f\"ur Mathematik, Gottfried-Claren-Str. 26,
53225 Bonn, Germany. $\quad$
e-mail: markan@mpim-bonn.mpg.de

\end{document}